\documentclass[review]{elsarticle}

\usepackage{lineno,hyperref}
\modulolinenumbers[5]
\usepackage[margin=1.0in]{geometry}
\usepackage{comment}
\usepackage{soul}
\usepackage{wrapfig}

\usepackage[cmex10]{amsmath}
\usepackage{amsfonts}
\usepackage{amssymb}
\usepackage{mathabx}
\usepackage{cases}
\usepackage{cancel}

\usepackage{caption}
\usepackage{subcaption}

\allowdisplaybreaks

\usepackage{amsthm}

\theoremstyle{remark}

\usepackage{amsmath}
\newcommand{\cT}{\mathcal{T}}

\newcommand{\cP}{\mathcal P}

\newcommand{\cK}{\mathbb {K}}

\usepackage{xcolor}
\usepackage{caption}
\usepackage{subcaption}

\usepackage{algorithm,algorithmic}

\usepackage{tabularx, booktabs}
\newcolumntype{L}{>{\raggedright\arraybackslash}X}
\newcolumntype{b}{X}
\newcolumntype{s}{>{\hsize=.5\hsize}X}
\usepackage{multirow}

\usepackage{enumitem}
\usepackage{nomencl}
\makenomenclature

\usepackage{etoolbox}
\renewcommand\nomgroup[1]{%
  \item[\bfseries
  \ifstrequal{#1}{B}{Parameters}{%
  \ifstrequal{#1}{C}{Variables}{%
  {%
  \ifstrequal{#1}{A}{Sets}{}}}}%
]}

\begin{document}

\begin{frontmatter}

\title{A Modeling Framework for Equitable Deployment of Energy Storage in Disadvantaged Communities}

\author{Miguel Heleno\fnref{lbnl}\corref{correspondingauthor}} 
\ead{miguelheleno@lbl.gov}
\author{Paul Lesur\fnref{us,CSGI}}
\author{Alexandre Moreira\fnref{lbnl}}

\address[lbnl]{Lawrence Berkeley National Laboratory, Berkeley, CA 94720, USA.}
\address[us]{Stanford University, Stanford, CA, USA.}
\address[CSGI]{Ecole Polytechnique, Palaiseau, France.}

\cortext[correspondingauthor]{Corresponding author}

\begin{abstract}
This paper provides an analytical framework to incorporate the deployment of behind-the-meter energy storage coupled with rooftop solar, and their associated revenue streams, in the context of equitable energy policy interventions. We propose an extension to the Justice40 optimization model by adding storage and incorporating more realistic solar compensation mechanisms, such as net-billing, which allows for temporal revenue differentiation and the economic viability of behind-the-meter energy storage devices. The extended model includes household-level PV plus storage co-deployment alongside existing interventions, such as weatherization, rooftop PV only, community solar, and community wind. From a modeling perspective, we propose a novel approximation method to represent storage operations and revenue streams without expanding the temporal dimension of model, thus maintaining its computational efficiency. The proposed model is validated using a case study in Wayne County, Michigan, involving 3,651 energy insecure households. 
\end{abstract}


\end{frontmatter}

\nomenclature[A, 1]{$\cT$}{Set of tracts, indexed by $\tau$}
\nomenclature[A, 2]{$\cP$}{Set of policy interventions. indexed generically by $p$ and referring to the following interventions \{w: weatherization, rts: rooftop solar, cs: community solar, cw: community wind, batt: batteries deployment\}.}

\nomenclature[A, 4]{$\cK_{p}^\tau$}{Set of household archetypes eligible for policy intervention $p$ in tract $\tau$, indexed by $k$.}
\nomenclature[A, 5]{$\cK_{w-f}$}{Set of household archetypes eligible for weatherization policies with $f$ as heating fuel type.}

\nomenclature[B, 2]{$Nb_k$}{Number of households of archetype $k$.}
\nomenclature[B, 3]{$\zeta_{\tau}$}{Solar productivity factor for tract $\tau$.}
\nomenclature[B, 4]{$\eta_{\tau}$}{Wind productivity factor for tract $\tau$.}
\nomenclature[B, 5]{$IC_{p}$}{Investment costs associated with policy intervention $p$.}
\nomenclature[B, 6]{$L_{p}$}{Lifetime of the investments associated with policy intervention $p$.}
\nomenclature[B, 7]{$Pel$}{Price of electricity.}
\nomenclature[B, 8]{$PV_{rem}$}{Photovoltaic remuneration}
\nomenclature[B, 9]{$WS_k$}{Percentage of weatherization savings in building $k$.}
\nomenclature[B, 10]{$E^{f}_k$}{Baseline expenditure of building $k$ associated with the consumption of fuel $f$.}
\nomenclature[B,11]{$\widebar{RTS_k}$}{ Maximum rooftop solar allowed in archetype building $k$.}
\nomenclature[B,12]{$\widebar{CS_\tau}$}{ Maximum community solar in tract $\tau$.}
\nomenclature[B,13]{$\widebar{CW_\tau}$}{ Maximum of community wind in tract $\tau$.}
\nomenclature[B,14]{$\widebar{B}$}{Maximum budget available.}
\nomenclature[B,15]{$Z_k^1$}{Rooftop solar capacity above which there is PV surplus for household $k$.}
\nomenclature[B,16]{$D^{batt}$}{Storage technology duration.}

\nomenclature[C, 1]{$d^w_{k}$}{Fraction of households of archetype $k$ to receive weatherization interventions.}
\nomenclature[C, 2]{$d^{rts}_{k}$}{Amount of rooftop PV installed in archetype  building $k$.}
\nomenclature[C, 3]{$d^{batt}_{k}$}{Amount of batteries installed in archetype building $k$.}
\nomenclature[C, 4]{$\delta^{batt}_{k}$}{Binary Variable describing whether batteries should be installed in archetype building $k$.}
\nomenclature[C, 5]{$d^{cs}_{t}$}{Amount of community PV installed in tract $t$.}
\nomenclature[C, 6]{$d^{cw}_{t}$}{Amount of community wind installed in tract $t$.}
\nomenclature[C, 7]{$g^{rts}_{k}$}{Amount of electricity generated from rooftop PV by an household archetype $k$.}
\nomenclature[C, 8]{$g^{cs}_{k}$}{Portion of electricity generated from community PV attributed to building $k$.}
\nomenclature[C, 9]{$g^{cw}_{k}$}{Portion of electricity generated from community wind attributed to building $k$.}
\nomenclature[C, 10]{$c^{rts}_{k}$}{Annualized costs of rooftop PV installations assigned to the archetype  household $k$.}
\nomenclature[C, 11]{$c^{batt}_{k}$}{Annualized costs of storage capacities installations assigned to the archetype  household $k$.}
\nomenclature[C, 12]{$c^{cs}_{t}$}{Annualized cost of community PV installation in tract $t$.}
\nomenclature[C,13]{$c^{cw}_{t}$}{Annualized cost of community wind installation in tract $t$.}
\nomenclature[C,14]{$eld_{k}$}{Electricity demand of an household archetype $k$.}
\nomenclature[C,15]{$ec_{k}$}{Total energy consumption of an household archetype $k$.}
\nomenclature[C,16]{$eg_{k}$}{Total energy generation associated with the household archetype $k$.}
\nomenclature[C,17]{$eb_{k}$}{Energy burden of an household archetype $k$.}
\nomenclature[C,18]{$\Delta eb^+_{k}$}{Energy overburden of an household archetype $k$, corresponding to the portion energy burden above the threshold.}
\nomenclature[C,19]{$sc^{rts}_k$}{Electricity produced by rooftop solar of household $k$ and self-consumed locally.}
\nomenclature[C,20]{$st^{rts}_k$}{Electricity produced by rooftop solar of household $k$ and stored in the battery for latter use.}
\nomenclature[C,21]{$sg^{rts}_k$}{Electricity produced by rooftop solar of household $k$ that is injected directly into the grid.}
\nomenclature[C,22]{$pd_{k,t},pc_{k,t}$}{Hourly charging/discarding to/from the battery installed at the household $k$.}
\nomenclature[C,23]{$soc_t^{k}$}{Energy state-of-charge of the battery installed at the household $k$.}
\nomenclature[C,24]{$z_k$}{Binary variable indicating whether there is surplus of rooftop solar.}

\printnomenclature

\section{Introduction}

\subsection{Motivation}

Energy costs have reached new records since the beginning of the year 2020 \cite{US_EIA}. This rapid increase has been mainly driven by oil and gas prices and impacted directly customers' purchasing power. As a consequence, energy insecurity among low-income households and restrictions to energy access \cite{memmott2021sociodemographic} have risen, leading to utility disconnections \cite{BAKER2021112663}. In the US, studies led by the American Council for Energy-Efficient Economy have demonstrated that low-income households as well as Black, Hispanic and Native American communities experience higher energy insecurity at the national, regional and metropolitan levels \cite{drehobl2020high}. While the threshold for a decent energy burden is commonly admitted to be at 6\%, low-income households with older adults and children can experience three-times higher energy burden. In several metropolitan areas such as Detroit and Boston, more than 25\% of low-income households have an energy burden reaching more than 15\%.

Geographic and socioeconomic disparities in energy burden in conjunction with the need to decarbonize the energy systems have motivated a new generation of place-based policy interventions \cite{osti_1900013, DOLTER2018221}. For example, the US Department of Energy launched the Local Energy Action Program (LEAP), which aims at addressing a clean and fast but equitable energy transition on the path towards sustainability \cite{osti_2000949}. Together with other programs, like the Weatherization Assistance program (WAP) or the Low-Income Home Energy Assistance Program (LIHEAP)\cite{brown2020high}, LEAP targets local communities with a specific share of low-income population to help them design policies to reduce energy consumption, greenhouse gas emissions and their energy burden. The benefits of these community level programs have been demonstrated and some improvements in the definition of eligible populations have been proposed \cite{murray2014impact}. 

As a strategy to alleviate the energy burden of low-income households, these programs often include several mechanisms for the deployment of renewable distributed generation, such as the installation of photovoltaic (PV) panels, together with energy efficiency improvements. In fact, the presence of rooftop solar generation devices behind-the-meter has the potential to significantly reduce household energy costs by increasing PV self-consumption, reducing overall electricity needs from the grid, and creating new revenue streams due to solar compensation mechanisms \cite{wan1996net}. 

Besides weatherization interventions and deployment of distributed generation, energy storage is also emerging as an efficient and cost-effective equitable investment that can be deployed at the residential level (behind-the-meter) or at the community level. In fact, batteries can help decrease energy cost, particularly in low-income households \cite{tarekegne2021energy}, and, when paired with rooftop PV, provide energy savings \cite{hoff2007maximizing} and backup power to consumers. They enhance households' electricity affordability and resilience \cite{MCNAMARA2022107063}, improve self-consumption \cite{lund2018capacity} and decrease overall household energy burden \cite{guan2023burden}. If provided with the right incentives, these PV plus storage installations can be even more impactful from an economic \cite{SANIHASSAN2017422}, reliability \cite{MAHESHWARI2020114964} and resilience \cite{GORMAN2023121166} perspective.

Given this potential of energy storage to benefit low-income communities, it is important that methodologies and tools to support the design of equitable energy policy interventions can properly capture those benefits and accurately inform the deployment of battery technologies in disadvantaged territories. Currently, there is no quantitative methodology that explicitly accounts for energy storage in equitable planning of place-based energy investments. This paper aims to address this issue, by proposing an extension of the Justice40 model, recently introduced in \cite{heleno2022optimizing}, to include PV plus storage installations in low-income communities.

\subsection{Literature Review}

In the field of equitable policy, a large body of literature has analyzed disparities across different sociodemographic groups regarding access to solar technologies \cite{Barbose2020, REAMES2020101612}, energy efficiency 
 \cite{FORRESTER2020114307}, weatherization interventions \cite{XU2019763, REAMES2016549} and energy storage \cite{BROWN2022112877}. These works, based on advanced empirical analyses, characterize inequities in energy infrastructure and emphasize the need to correct them through policy interventions. However, these empirical analyses cannot be used to design forward-looking policies nor to plan the deployment of energy assets. For that purpose, it is necessary to use techno-economic models that can simulate future scenarios and support decisions around energy infrastructure. Those models are part of another separate body of literature, focused on forward-looking energy resource planning in a variety of applications, including integrated energy systems \cite{VANBEUZEKOM2021116880}, weatherization interventions \cite{ROGEAU2020114639} or microgrids \cite{dercam} \cite{Homer}. Nevertheless, these models are exclusively focused on techno-economic aspects and they fail to capture the sociodemographic dimension and, therefore, social inequities in energy infrastructure deployment.

The gap between these two bodies of literature was addressed by a techno-economic modeling framework for the design of equitable energy interventions that explicitly captures the sociodemographic aspects \cite{heleno2022optimizing}. This framework comprises the methodology behind the open-source Justice40 tool \cite{Justice40Tool}, developed by the Lawrence Berkeley National Laboratory for the US Department of Energy, and used to support decisions around deployment of distributed generation and weatherization interventions. However, this model does not include the deployment of battery technologies as part of the portfolio of energy resources, neglecting the social benefits of PV plus storage installations discussed above \cite{tarekegne2021energy}. This limitation is particularly relevant as compensation mechanisms for residential solar are evolving toward net-billing models that incentivize using storage to maximize solar self-consumption \cite{FORRESTER2022104714}.
 
The technical reason for this limitation is related to the characteristics of the energy storage models, which usually track energy generation and consumption over time to compute the state of charge and the charging/discharging operations of battery devices \cite{sioshansi2021energy}. These operation models can then be extended to implement new features or new objective functions. For example, in \cite{barbour2018projecting} and \cite{hassan2017optimal}, the authors present a model that provides the optimal energy storage sizing and management that minimize the overall cost of electricity for consumers. To do so, the authors use time series of households electricity consumption, PV generation and electricity prices. This allows to explicitly model temporal price-differentiation and capture storage-specific revenue streams, such as arbitrage or self-consumption. 

These temporal representation of storage devices implies capturing operations at the intra-day level (typically 1 hour or 15 minutes), which, in optimization models, requires thousands of new variables and constraints to represent a single storage device per year. Therefore, such extension can significantly increase the computational time in applications that consider a large number of potential storage installations. This is the case of energy planning models with equity objectives that entail a very granular representation of sociodemographic groups and their energy uses \cite{heleno2022optimizing}. So far, the only feasible way of combining storage operations with granular sociodemographic representations is to run models offline, during a long time, and then include their results in policy analyses \cite{GORMAN2023121166}. 
\subsection{Contributions}

This work proposes an extension of the justice40 model \cite{heleno2022optimizing} and its online application \cite{Justice40Tool} to account for PV plus storage and their revenue streams in the context of equitable deployment of energy policy interventions. The contributions are threefold:
\begin{itemize}
    \item We extend the energy burden formulation to consider more realistic solar compensation mechanisms that imply temporal revenue differentiation, such as net-billing, and create an economic opportunity for behind-the-meter energy storage devices.
    \item We add household-level PV plus storage co-deployment as a part of the portfolio of equitable interventions represented in the Justice40 model, which includes weatherization, rooftop PV only, as well as community solar and community wind. 
    \item We propose an approximation of PV plus storage representation that allows to capture their operations and revenue streams without using a temporal expansion of the model and, therefore, keeping the computational burden compatible with online applications.
\end{itemize}

To illustrate and validate the the proposed model, we compare it with a version comprising a full temporal representation of storage operations, using a case study in Wayne County, Michigan, US, involving a population of 3,651 energy insecure households.

\subsection{Organization of the paper}

The rest of the paper is divided as follows: section \ref{sec:methodology} presents an overview of the existing justice40 model. Section \ref{sec:new_m} introduces the re-formulation and the temporal approximations proposed to represent solar plus storage in the model. Section \ref{sec:results} compares the output of the proposed model with a complete time-dependent formulation and discussed the results.

\section{Base optimization framework: the Justice40 model}
\label{sec:methodology}

\subsection{Model overview}

As stated above, this paper extends the optimization framework presented in \cite{heleno2022optimizing}, developed under the US Department of Energy Justice 40 initiative \cite{Justice40}. The Justice40 model aims at determining the optimal portfolio of energy interventions to reduce the energy insecurity in a particular geographical location. Currently, the model includes four types of energy interventions: weatherization of residential buildings, deployment of rooftop solar panels, and deployment of community-owned solar and wind. The model is currently open to the public both through an interface and an open-source python module \cite{Justice40Tool}.  

Energy interventions in the model are deployed in an equitable manner with the objective of minimizing the energy insecurity of an eligible population, subjected to a budget constraint ($\Bar{B}$). Within this framework, the energy insecurity of a particular household is characterized by an excessive energy burden above a certain threshold ($\Delta eb^+_{k}$), often considered to be approximately 6\%. Formally, the objective function and the budget constraints are represented by equations (\ref{eqn:objective_alternative}) and (\ref{eqn:budget constraints}), respectively. The population is represented by a set of households archetypes ($\cK$), which capture different households conditions across a set of U.S Census tracts ($\cT$). 

 \begin{align}
& min ~ \sum_{k \in \cK} \Delta eb^+_{k} \cdot Nb_k \label{eqn:objective_alternative} \\
& \sum_{k \in \cK_{w}} c_k^{w} \cdot Nb_k +  \sum_{k \in \cK_{rts}} c^{rts}_k \cdot Nb_k + \sum_{\tau \in \cT} (c^{cs}_\tau + c^{cw}_\tau) - \theta \Bar{B} \leq 0
\label{eqn:budget constraints}
\end{align}

The annualized costs of deploying household-level interventions ($hl$) - rooftop solar ($rts$) and weatherization ($w$) - and community-owned, tract-level, interventions ($cl$) - solar ($cs$) and wind ($cw$) - are given by equations (\ref{eqn:household_level}) and (\ref{eqn:comlevel}), respectively. The distributed generation allocated to each household archetype is expressed by constraints (\ref{eqn:rtop_generation})-(\ref{eqn:energy_generation}) as a function of solar ($\zeta_{\tau}$) and wind ($\eta_{\tau}$) capacity factors in each tract. Equations (\ref{eqn:energy_costs_weatherization_nel}) express the electricity costs for buildings with non-electric heating-fuel ($hf$), whose consumption is affected by the deployment of weatherization ($d^{w}_k$) and corresponding savings ($WS_k$). For those buildings, the electricity demand remains constant, as captured in constraints (\ref{eqn:consumption_weatherization_nel}). Equations (\ref{eqn:consumption_weatherization_ele}) model the opposite situation, in which electricity is the heating source and the electric consumption is affected by the weatherization interventions. In those cases, the costs of non-heating, non-electric, fuels ($of$) remains constant, as expressed in constraints (\ref{eqn:energy_costs_weatherization_ele}).

\begin{align}
& c^{hl}_k = IC_{hl} \cdot d^{hl}_k \cdot \frac{r}{1-(1+r)^{-L_{hl}}} \quad \forall k \in \cK_{hl} 
\label{eqn:household_level} \\
& c^{cl}_\tau = IC_{cl}  \cdot d^{cl}_\tau  \cdot \frac{r}{1-(1+r)^{-L_{cl}}} \quad \forall \tau \in \cT \label{eqn:comlevel}\\
& g^{rts}_k =d^{rts}_k \cdot \zeta_{\tau} \quad \forall k \in \cK_{rts}
\label{eqn:rtop_generation} \\
& d^{rts}_k \leq \overline{RTS_k} \quad \forall k \in \cK_{rts}
\label{eqn:rtop_generation?lim} \\
& g_k^{cs} = \frac{d^{cs}_{\tau}\cdot \zeta_{\tau}}{\sum_{k \in \cK_{cs}^\tau} Nb_k} \quad \forall k \in \cK_{cs}^\tau \quad \forall \tau \in \cT 
\label{eqn:comsolar_generation} \\
& g_k^{cw} = \frac{d^{cw}_{\tau}\cdot \eta_{\tau}}{\sum_{k \in \cK_{cw}^\tau} Nb_k} \quad \forall k \in \cK_{cw}^\tau \quad \forall \tau \in \cT \label{eqn:community_wind_geneartion}\\
& eg_{k} =g^{rts}_k + g^{cs}_k + g^{cw}_k \quad \forall k \in \cK \label{eqn:energy_generation}\\
& ec_{k} = eld_{k} \cdot Pel + E^{hf}_k - (E^{hf}_k \cdot  d^{w}_k \cdot WS_k) + E^{of}_k \quad \forall k \in \cK_w \setminus \cK^{el}_w  \label{eqn:energy_costs_weatherization_nel} \\
& eld_{k} = \frac{E^{el}_k}{Pel}  \quad \forall k \in \cK_w \setminus \cK^{el}_w 
\label{eqn:consumption_weatherization_nel} \\
& eld_{k} = \frac{E^{el}_k}{Pel} -  d^{w}_k \cdot WS_k \cdot \frac{E^{el}_k}{Pel} \quad \forall k \in \cK^{el}_w 
\label{eqn:consumption_weatherization_ele} \\
& ec_{k} = eld_{k} \cdot Pel + E^{of}_k \quad \forall k \in \cK^{el}_w \label{eqn:energy_costs_weatherization_ele}
\end{align} 

Finally, equations (\ref{eqn:energy burden_definition}) model energy burden of each household archetype, given by the net energy expenditures (i.e., the difference between energy costs and distributed generation revenues) divided by the annual income of the household $I_k$. In equation (\ref{eqn:burden_threshold}), the energy burden is decomposed into positive and negative deviations in relation to a threshold ($\Bar{Eb}$), above which a household is considered energy insecure. Therefore, addressing energy insecurity is equivalent to minimizing the positive deviations regarding the energy insecurity threshold, as presented in the objective function (\ref{eqn:objective_alternative}).

\begin{align}
& eb_{k} = \frac{ec_k - eg_{k} \cdot Pel}{I_k} \quad \forall k \in \cK \label{eqn:energy burden_definition} \\
& eb_{k} - \Bar{Eb} = \Delta eb^+_{k} - \Delta eb^-_{k} \quad \forall k \in \cK \label{eqn:burden_threshold} \\
& \Delta eb^+_{k}, \Delta eb^-_{k}, eb_{k} \geq 0 \quad \forall k \in \cK \label{eqn:burden_delta_lim}
\end{align} 

Here, we presented a short summary of the justice40 model in its current state, using the variables and constraints that are relevant for the expansion proposed in this paper. The detailed formulation of the model is discussed in \cite{heleno2022optimizing}.

\subsection{Limitations to storage deployment modeling}

As formulated, this model is unable to properly capture the equitable benefits related to the deployment of energy storage in combination with rooftop solar. In fact, the model lacks temporal representation and, in the definition of energy burden (\ref{eqn:energy burden_definition}), implicitly assumes a net-metering revenue associated with distributed generation. These characteristics do not allow the introduction of more realistic revenue functions, such as net-billing, that include temporal price differentiation and, therefore, create conditions for consumers to benefit from energy storage.

\section{Model re-formulation to include storage}
\label{sec:new_m}

\subsection{Rooftop solar and storage co-deployment}

In  the new version of the model, households eligible for the deployment of rooftop solar may also be eligible for the deployment of behind-the-meter battery installations. Often these PV plus storage solutions are offered by the industry in standardized sizes, comprising a battery packages proportional to capacity of the PV panel. Therefore, in our co-deployment, a linear relation between
the batteries and the rooftop solar panel sizes is assumed. The parameter $\beta$, in equation (\ref{eqn:linear_rel_PV_batt}) translates this relationship and allows to write a potential battery deployment as a function of the rooftop PV.  

\begin{equation}
\beta =\frac{d_{k}^{batt}}{d_{k}^{rts}} \quad \forall k \in \cK_{rts}
\label{eqn:linear_rel_PV_batt}
\end{equation}

The annualized cost of battery installation can also express the linear relationship of co-deployment with PV, as shown in equation (\ref{eqn:cost_battery}), written as an equivalent of constraints (\ref{eqn:household_level}) and (\ref{eqn:comlevel}) to storage technologies. It is important to note that the battery installation remains as an investment option, captured by the binary variable $\delta_k^{batt}$, to allow the possibility of not considering storage at all. 

\begin{align}
& c^{batt}_k = {IC}^{batt}_{k} \cdot \delta_k^{batt} \cdot \beta \cdot d_k^{rts} \cdot \frac{r}{1-(1+r)^{-L_{batt}}} \label{eqn:cost_battery}  
\end{align}

As written above, the cost function implies the multiplication of two variables: the binary decision associated with the presence of batteries ($\delta_k^{batt}$) and the solar deployment deployment ($d_k^{rts}$). This nonlinear relationship can be easily addressed by explicitly modeling the battery capacity in a disjunctive manner and bounding it using the rooftop solar limit at the household level ($\widebar{RTS_k}$). Equations (\ref{eqn:cost_batteries_linear})-(\ref{eqn:cost_batteries_sup}) rewrite the battery cost and capacity constraints in an equivalent form.

\begin{align}
& c^{batt}_k = {IC}^{batt}_{k} \cdot cap_k^{batt} \cdot \frac{r}{1-(1+r)^{-L_{batt}}} \label{eqn:cost_batteries_linear} \\
& cap_k^{batt} \leq \beta \cdot d_k^{rts}  \\
& cap_k^{batt} \geq \beta \cdot d_k^{rts} - \widebar{RTS_k} \cdot \beta \cdot (1-\delta_k^{batt}) \label{eqn:cost_batteries_dij_inf} \\
& cap_k^{batt} \leq \widebar{RTS_k} \cdot \beta \cdot \delta_k^{batt} \label{eqn:cost_batteries_sup} \\
& cap_k^{batt} \geq 0 \label{eqn:cost_batteries_inf}
\end{align}

\subsection{Energy burden with different solar compensation mechanisms}

The original model implicitly assumes a net-metering scheme, where PV sales to the grid are compensated at the same rate as electricity demand. However, new solar compensation schemes, such as net-billing, pay a lower rate for PV electricity than the demand rate. This creates an economic opportunity for behind-the-meter batteries, which can store excess PV energy instead of selling it back to the grid at a lower price, allowing it to be used later to offset demand.

To model the price differentiation in solar compensation, equation (23) decomposes rooftop solar generation into three components: the first corresponds to the electricity produced by the system and self-consumed locally ($sc^{rts}_k$) to supply the electric demand; the second relates to the PV generation stored in the battery for later use ($st^{rts}_k$); and the third corresponds to the PV generation injected into the grid ($sg^{rts}_k$).

\begin{equation}
g^{rts}_k =sc^{rts}_k + st^{rts}_k + sg^{rts}_k \quad \forall k \in \cK_{rts}
\label{eqn:rooftop_generation}
\end{equation} 

Within a net-billing scheme, the self-consumed portion of solar generation offsets the demand, effectively remunerating it at the electricity price ($Pel$). The same applies to PV generation stored in the battery and later used to offset demand. In contrast, the remaining solar generation is sold to the grid at a specific PV remuneration rate ($PVrem$). This rate defines the compensation scheme: net-billing occurs when the PV remuneration is less than the electricity price ($PVrem<Pel$), while net-metering occurs when they are equal ($PVrem=Pel$).

The total rooftop solar compensation, including the effect of the battery, is used to rewrite the energy burden equation (\ref{eqn:energy burden_new}), considering different PV remunerations for generation consumed locally and injected into the grid. To ensure consistency with the previous version of the model, generation from community solar ($g^{cs}_k$) and community wind ($g^{cw}_k$) remains remunerated based on the net-metering scheme. Depending on the regulatory framework, different remunerations could be introduced for these two types of generation. However, this extension is straightforward and is not covered in this paper.

\begin{align}
& eb_{k} = \frac{ec_k - (sc^{rts}_k + st^{rts}_k + g^{cs}_k + g^{cw}_k) \cdot Pel - sg^{rts}_k \cdot PVrem}{I_k} \quad \forall k \in \cK \label{eqn:energy burden_new}
\end{align}

\subsection{Rooftop solar and storage operation}

Expanding the formulation to include different remuneration schemes and storage investments requires accurate modeling of the three components of solar generation: self-consumption ($sc^{rts}_k$), storage ($st^{rts}_k$), and grid injection ($sg^{rts}_k$). These variables depend not only on the amount of solar and storage deployment but also on operational aspects over time. Self-consumption ($sc^{rts}_k$) relies on the temporal correlation between solar production and electricity demand. During hours of surplus PV production, battery operation ($st^{rts}_k$) can shift energy to later times, helping to decrease the household's energy burden, as shown in equation (\ref{eqn:energy burden_new}). When the battery capacity is insufficient to store all the surplus PV energy, the remaining generation is sold to the grid.

\subsubsection{Time-resolved model}

A direct approach to capturing these temporal components of rooftop generation is to expand the variables in the time domain ($T$) and explicitly model the battery operation. Assuming that $pd_{k,t}$ is the hourly discharge of the battery and $sc^{rts}_{k,t}$, $sg^{rts}_{k,t}$ are expansions in $T$ of the self-consumption and PV sales, respectively, the three rooftop generation components can be given by equations (\ref{eq:stagg})-(\ref{eq:sgagg}). The equations (\ref{eq:solar_time}) expand in the time domain the solar generation condition presented in (\ref{eqn:rooftop_generation}), using the battery charging ($pc_{k,t}$). Equations (\ref{eq:sc_time}) impose the definition of self-consumption as a function of the electricity demand at each time ($eld_{k,t}$). Equations (\ref{eq:soc}) model the battery state of charge. Finally, the constrains (\ref{eq:lim_ch}) and (\ref{eq:lim_dc}) impose the charging and discharging limits at each time.

\begin{align}
& st^{rts}_k = \sum_{t \in T} pd_{k,t} \quad \forall k \in \cK_{batt} \label{eq:stagg}\\
& sc^{rts}_k = \sum_{t \in T} sc^{rts}_{k,t} \quad \forall k \in \cK_{batt} \label{eq:scagg}\\
& sg^{rts}_k = \sum_{t \in T} sg^{rts}_{k,t} \quad \forall k \in \cK_{batt} \label{eq:sgagg}\\
& sc^{rts}_{k,t} + sg^{rts}_{k,t} + pc_{k,t} =  d^{rts}_k \cdot \zeta_{\tau,t} \quad \forall k \in \cK_{batt} \label{eq:solar_time}\\ 
& sc^{rts}_{k,t} \leq eld_{k,t} \quad \forall k \in \cK_{batt} \label{eq:sc_time} \\ 
& soc_{k,t} = soc_{k,t} + pc_{k,t} -pd_{k,t}  \quad \forall k \in \cK_{batt} \label{eq:soc}  \\
& soc_{t,k}  \le cap_k^{batt} \quad \forall k \in \cK_{batt} \label{eq:soc_lim} \\
&  0 \le  pc_{k,t} \cdot D^{batt} \le  cap_k^{batt} \quad \forall k \in \cK_{batt} \label{eq:lim_ch}\\
&  0 \le  pd_{k,t} \cdot D^{batt} \le  cap_k^{batt} \quad \forall k \in \cK_{batt} \label{eq:lim_dc}
\end{align}

Although the temporal expansion can accurately capture solar plus storage operation, it significantly increases the number of variables in the model. Since the model is now a Mixed Integer Linear Program (MILP) due to the introduction of the binary battery investment option ($\delta_k^{batt}$), this increase in variables introduces new computational burdens. It is important to note that the model resolution is at the level of household archetypes in each tract, which implies thousands of archetypes for most U.S. counties. 
Therefore, a direct expansion in the temporal domain would exponentially increase the computational burden, making it incompatible with a user-friendly application such as the Justice40 tool \cite{Justice40Tool}.

\subsubsection{Proposed self-consumption approximation}

An alternative to a full representation of the temporal solar plus storage operation is to describe the three components of annual solar generation — self-consumed ($sc^{rts}_k$), stored ($st^{rts}_k$), and grid-injected ($sg^{rts}_k$) — as a function of rooftop solar deployment, as shown in (\ref{eqn:components}). In fact, such function is possible if we assume a fixed PV plus storage size relationship ($\beta$) and a fixed battery duration ($D$). It is important to stress that these functions do not require directly representing the temporal operations of the solar and battery resources. However, the functions associated with the self-consumed, $sc^{rts}_k (d_{k}^{rts})$, and stored, $st^{rts}_k (d_{k}^{rts})$, components of solar generation need to be obtained.

\begin{equation}
     \begin{cases}
        sc^{rts}_k (d_{k}^{rts}) =   f^{sc}_k(d_{k}^{rts}) \\
        st^{rts}_k (d_{k}^{rts}) =   f^{st}_k(d_{k}^{rts}) \\
        sg^{rts}_k (d_{k}^{rts}) = d^{rts}_k \cdot \zeta_{\tau} -   st^{rts}_k (d_{k}^{rts}) -  sc^{rts}_k (d_{k}^{rts}) 
     \end{cases}
\label{eqn:components}
\end{equation}

\begin{figure}[h]
         \centering
         \includegraphics[scale = 0.3]{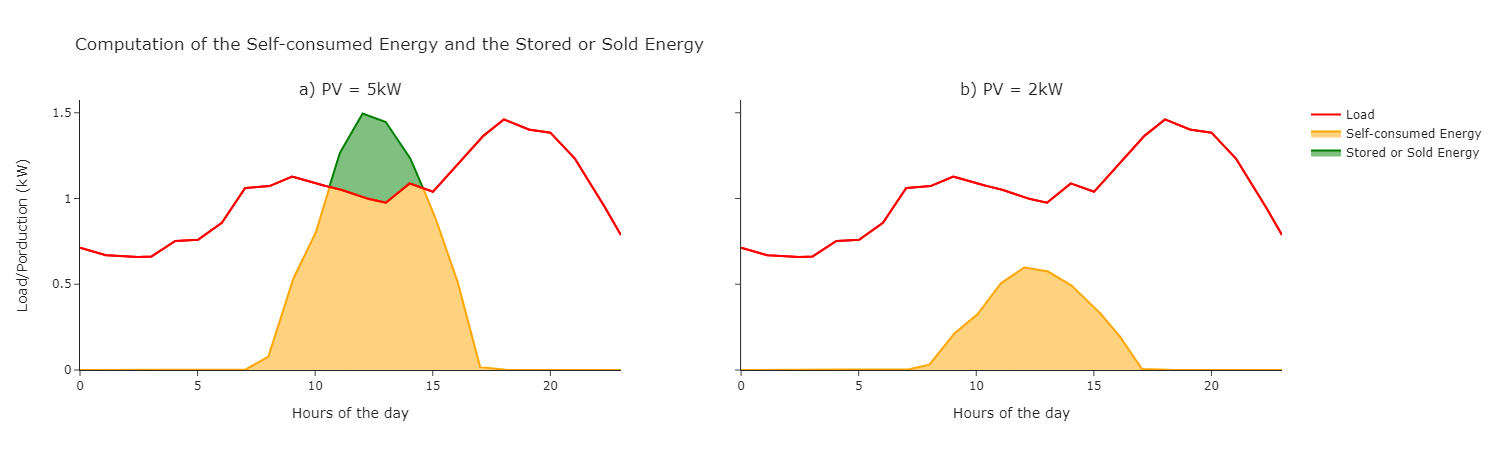}
        \caption{Self-consumed solar generation for two cases of rooftop solar deployment: a) when $d_{k}^{rts} \geq Z_k^1$; b) when $d_{k}^{rts} < Z_k^1$.}
        \label{figure:solar_sc_and_grid}
\end{figure}

To describe the energy self-consumed as a function of the rooftop deployment, $sc^{rts}_k (d_{k}^{rts})$, we take the yellow area of solar generation under the demand curve illustrated in Figure \ref{figure:solar_sc_and_grid}. As shown in the left-hand side of the figure, there is a quantity of solar deployment ($d_{k}^{rts}=Z_k^1$) above which a surplus generation occurs in, at least, one period of the year (the figure only shows one day of operation for simplicity). When the rooftop deployment is below $Z_k^1$ (right-hand side of the figure), the entire PV generation is self-consumed. Hence, the PV self-consumption can be described as piecewise function of the rooftop deployment, illustrated by equation (\ref{eqn:rts_sc}), accounting for capacity values lower and higher than $Z_k^1$.

\begin{equation}
sc^{rts}_k (d_{k}^{rts}) =
     \begin{cases}
        d^{rts}_k \cdot \zeta_{\tau} & \text{if } d_{k}^{rts} \leq Z_k^1\\
        f^{2,sc}_k(d_{k}^{rts}) & \text{if } d_{k}^{rts} > Z_k^1
     \end{cases}\quad \forall k \in \cK_{rts}
\label{eqn:rts_sc}
\end{equation}

\begin{figure}[h]
         \centering
         \includegraphics[scale = 0.6]{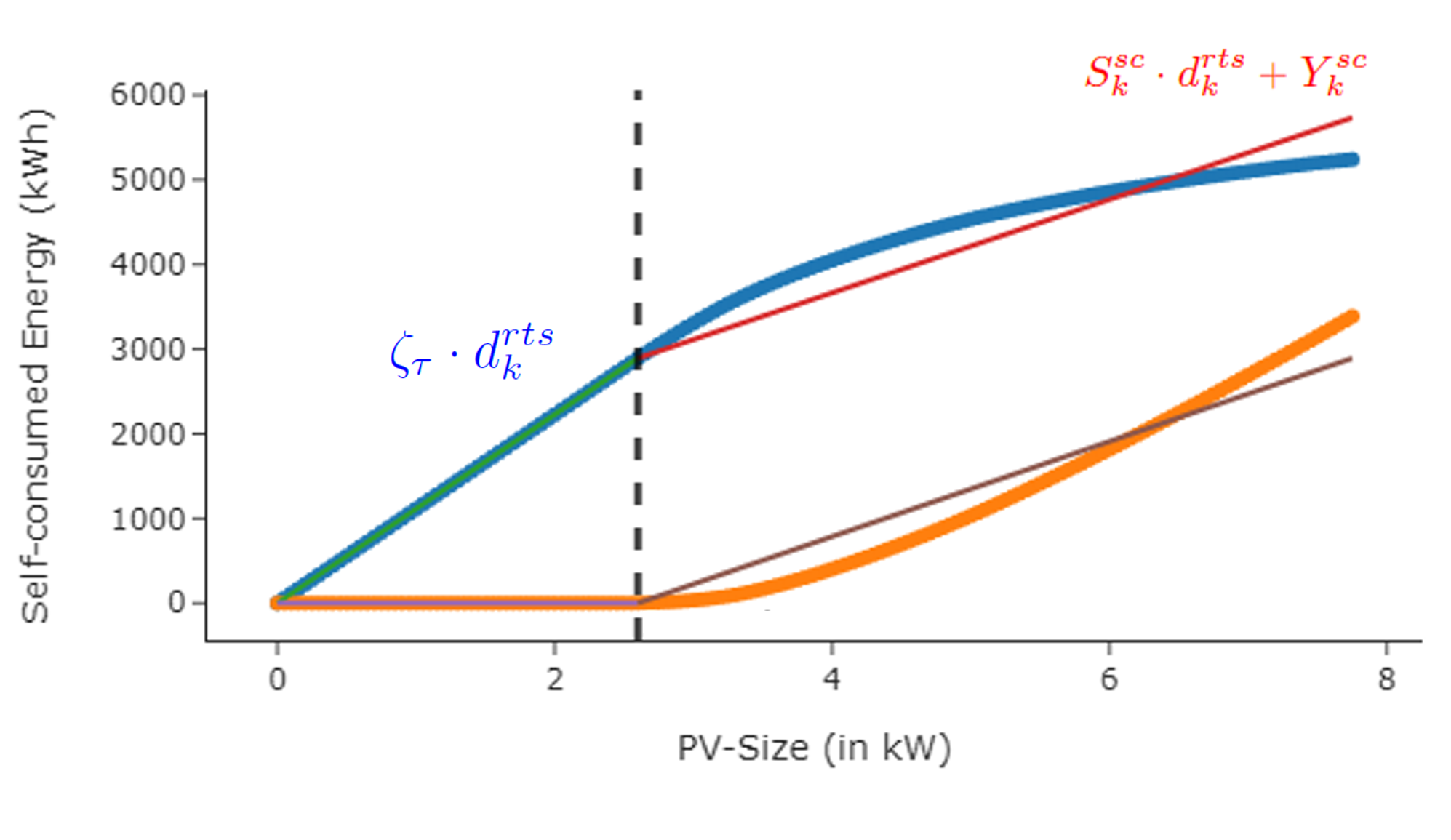}
        \caption{Self-consumption as a function of PV deployment, $sc^{rts}_k (d_{k}^{rts})$ }
        \label{figure:solar_sc_linearization}
\end{figure}

Although the first piece of the function (\ref{eqn:rts_sc}) scales linearly with the capacity factor of each tract ($\zeta_{\tau}$), the second piece ($f^{2,sc}_k$) is more complex to obtain, as it entails describing the yellow area under the demand curve shown on the left-hand side of the Figure \ref{figure:solar_sc_and_grid}. However, given a household ($k$) with an electricity demand profile and a value of $Z_k^1$, the self-consumption area can be pre-computed  for different levels of rooftop PV deployment ($d_{k}^{rts}$), sampled within the interval $[0,\widebar{RTS_k}]$. The blue curve in Figure \ref{figure:solar_sc_linearization} depicts the full function $sc^{rts}_k (d_{k}^{rts})$ with the two pieces described above. While the first piece is linear, $f^{2,sc}_k$ is not, but it can be approximated by a linear regression, described by a slope ($S^{sc}_k$) and an intercept ($Y^{sc}_k $). 

Considering this linear approximation ($f^{2,sc}_k \approx S^{sc}_k \cdot d_{k}^{rts} + Y^{sc}_k$), the two components of the self-consumption function, $sc^{rts}_k (d_{k}^{rts})$, can be implemented in a disjunctive manner, considering the value of $Z_k^1$, as described in equations (\ref{zg})-(\ref{1l}). Constraints (\ref{zg})-(\ref{d_bound}) impose the limits of the deployment according to the value of $Z_k^1$. In other words, if the binary variable $z_k$ is equal to 1, the second piece of $sc^{rts}_k (d_{k}^{rts})$ is activated, which implies $d_{k}^{rts} (\ref{zg})-(\ref{1l}) \geq Z_k^1$, while for $z_k=0$, the rooftop deployment is constrained by $Z_k^1$. Constraints (\ref{2g}) and (\ref{2l}) ensure consistency with the value of $Z_k^1$ and impose the linear approximation of self-consumption when $z_k=1$. Analogously, constraints (\ref{1g}) and (\ref{1l}) impose first piece of the self-consumption function described in (\ref{eqn:rts_sc}). 

\begin{align}
& d_{k}^{rts} \geq Z_k^1 - \widebar{RTS_k} \cdot (1-z_k) \label{zg}\\
& d_{k}^{rts} \leq Z_k^1 + \widebar{RTS_k} \cdot z_k \label{zl} \\ 
& d_{k}^{rts} \geq 0 \label{d_bound} \\
& sc^{rts}_k \geq S^{sc}_k \cdot d_{k}^{rts} + Y^{sc}_k - \zeta_{\tau} \cdot \widebar{RTS_k} \cdot T \cdot (1-z_k) \label{2g}\\
& sc^{rts}_k \leq S^{sc}_k \cdot d_{k}^{rts} + Y^{sc}_k + \zeta_{\tau} \cdot \widebar{RTS_k} \cdot T \cdot (1-z_k)\label{2l}\\
& sc^{rts}_k \geq \zeta_{\tau} \cdot d_{k}^{rts} - \zeta_{\tau} \cdot \widebar{RTS_k} \cdot T \cdot z_k \label{1g}\\
& sc^{rts}_k \leq \zeta_{\tau} \cdot d_{k}^{rts} + \zeta_{\tau} \cdot \widebar{RTS_k} \cdot T \cdot z_k  \label{1l}
\end{align}

\subsubsection{Energy Storage}

In addition to self-consumption, the load and PV profiles can also be used to derive the energy stored in a battery, $st^{rts}_k (d_{k}^{rts})$, using a similar regression approach. However, since the energy stored also depends on the battery deployment decision ($\delta_{k}^{batt}$), we use $\widebar{st^{rts}_k}$ to denote the opportunity for energy storage that can be materialized if a battery exists. Similarly to the self-consumption, this component of the solar generation can be written as a piecewise function of the rooftop deployment ($d_{rts}^k$) as demonstrated in equation (\ref{eqn:g_batt}). In fact, under the PV capacity $Z_k^1$, there is no generation surplus to be stored as all the electricity production from the rooftop solar is directly consumed by the household. When PV deployment is above $Z_k^1$, the opportunity to store generation surplus is given by $\widebar{st^{rts}_k}(\beta,d_{rts}^k)$.

\begin{equation}
\widebar{st^{rts}_k} (d_{k}^{rts}) =
     \begin{cases}
        0 & \text{if } d_{k}^{rts} < Z_k^1\\
        \widebar{st^{rts}_k}(\beta,d_{rts}^k) & \text{if } d_{k}^{rts} \geq Z_k^1
     \end{cases}\quad \forall k \in \cK_{rts}
\label{eqn:g_batt}
\end{equation}

\begin{figure}[h]
         \centering
         \includegraphics[scale = 0.35]{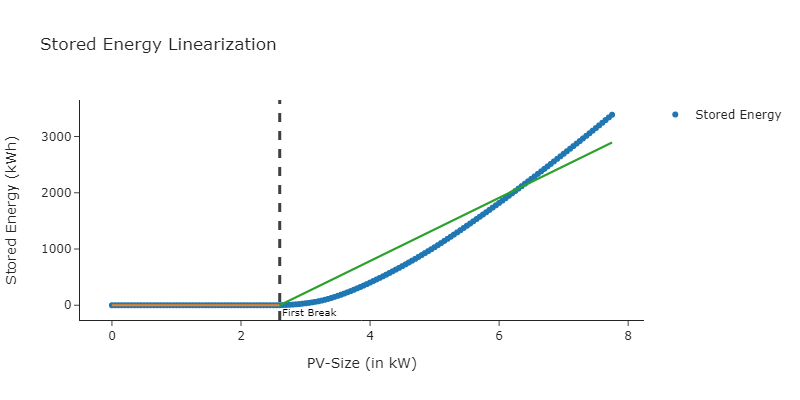}
        \caption{Linearization of the energy stored for one household}
        \label{fig:g_batt}
\end{figure}

As done for the self-consumption, it is possible calculate the opportunity for energy storage $\widebar{st^{rts}_k}(\beta,d_{rts}^k)$, for each household ($k$) and battery technology duration ($D^{batt}$). This entails pre-computing the battery dispatch (\ref{eq:stagg})-(\ref{eq:sgagg}) for multiple levels of rooftop deployment. Given the characteristics of the objective function, such dispatch is equivalent to a simple rule: the battery should charge as much as battery duration ($D^{batt}$) allows using the PV surplus and discharge later in the day to offset the demand. By pre-computing this dispatch for multiple values of rooftop deployment ($d_{rts}^k$) and corresponding battery capacities ($\beta \cdot d_{rts}^k$), we can obtain the function that characterizes the opportunity for energy storage,  $\widebar{st^{rts}_k}$, as illustrated in Figure \ref{fig:g_batt}. Similarly to self-consumption energy, it is possible to observe that this function $\widebar{st^{rts}_k}(\beta,d_{rts}^k)$ is not linear, but it can be approximated by a linear regression, described by a slope ($S^{st}_k$) and an intercept ($Y^{st}_k $).

This linearized quantity of solar generation stored in the battery, $st^{rts}_k$, is imposed by the constraints (\ref{lim_delta})-(\ref{st_lim}) in a disjunctive manner, using the battery deployment decision ($\delta_{k}^{batt}$) and the variable $z_k$, presented above, which models whether the solar deployment is above or below $Z_k^1$. 

\begin{align}
& \delta_{k}^{batt} \leq z_k \label{lim_delta}\\ 
& \widebar{st^{rts}_k} \geq S^{st}_k \cdot d_{k}^{rts} + Y^{st}_k - \zeta_{\tau} \cdot \widebar{RTS_k} \cdot T \cdot (1-z_k) \label{2_bat_g}\\
& \widebar{st^{rts}_k} \leq S^{st}_k \cdot d_{k}^{rts} + Y^{st}_k + \zeta_{\tau} \cdot \widebar{RTS_k} \cdot T \cdot (1-z_k)\label{2_bat_l}\\
& \widebar{st^{rts}_k}  \geq - \zeta_{\tau} \cdot \widebar{RTS_k} \cdot T \cdot z_k \label{1_bat_g}\\
& \widebar{st^{rts}_k} \leq \zeta_{\tau} \cdot \widebar{RTS_k} \cdot T \cdot z_k \label{1_bat_l}\\
& st^{rts}_k  \geq \widebar{st^{rts}_k} - \zeta_{\tau} \cdot \widebar{RTS_k} \cdot (1-\delta_{k}^{batt}) \label{1_st_g}\\
& st^{rts}_k \leq \widebar{RTS_k} \cdot \delta_{k}^{batt} \label{1_st_l}\\
& st^{rts}_k \leq \widebar{st^{rts}_k} \label{st_opp}\\
& \widebar{st^{rts}_k}, st^{rts}_k \geq 0 \label{st_lim} 
\end{align}

Constraints (\ref{lim_delta}) capture the relationship between these two variables: for solar deployment below $Z_k^1$ ($z_k$=0), there is no surplus and, therefore, storage can be neglected; for rooftop solar deployment above $Z_k^1$ ($z_k$=1), the battery deployment becomes optional and, therefore, the decision variable ($\delta_{k}^{batt}$) is allowed to assume a value equal to one. Constrains (\ref{2_bat_g}) and (\ref{2_bat_l}) ensure consistency in this this latter case, and impose the linear regression for storage opportunity when there is surplus in the system. In contrast, constrains (\ref{1_bat_g}) and (\ref{1_bat_l}) model the former case, when no surplus exists, forcing the opportunity for storage to be 0. Finally, constraints (\ref{1_st_g})- (\ref{st_lim}) model the actual stored energy as a function of the deployment decision ($\delta_{k}^{batt}$): when $\delta_{k}^{batt}=0$, no energy can be stored regardless of the storage opportunity value; when $\delta_{k}^{batt}=1$ the energy storage is limited by the storage opportunity ($\widebar{st^{rts}_k}$).

\section{Results}
\label{sec:results} 

In this section we apply the proposed methodology and the consequent enhancement of the Justice40 model to a population of energy insecure households. The results exclusively cover the contributions of this paper and are focused on the accuracy of the model and on the equity role of energy storage co-deployed with rooftop solar. General results on place-based energy equity investments and the role of other technologies are presented in the original Justice40 model paper \cite{heleno2022optimizing}. 

\subsection{Datasets}
The results presented in this section are based on model inputs sourced from public datasets. Non-temporal inputs are derived from the same datasets as described in the original Justice40 model \cite{heleno2022optimizing}. To support policy decisions regarding equitable interventions, the model incorporates inputs such as household archetypes and their locations using the Low-income Energy Affordability Data Tool (LEAD) \cite{Lead_Data}. Weatherization costs and savings projections rely on evaluations from the Weatherization Assistance Program \cite{WAP_eval}. Solar capacity factors are obtained from the Rooftop Energy Potential of Low Income Communities in America (REPLICA) \cite{REPLICA}, while wind capacity factors are derived from a study assessing onshore wind resource potential across the contiguous US \cite{LOPEZ2021120044}. Further details on how these datasets are utilized can be found in \cite{heleno2022optimizing}.

To construct the approximations of solar plus storage operation proposed in this paper, hourly generation and load profiles at the household archetype level are key inputs. Solar generation profiles include data from the NREL National Solar Radiation Database (NSRDB) \cite{sengupta2018national}, with PV productivity profiles generated using the System Advisory Model (SAM) \cite{blair2014system}. On the demand side, household electricity profiles are derived from ResStock datasets specific to each location \cite{oedi_4520}. For consistency with other datasets, solar generation and load profiles are adjusted to align with the annual solar capacity factor from REPLICA and electricity consumption data from the LEAD tool.

\subsection{Case Study}

We use the model to optimize equitable policy interventions, including energy storage, in Wayne County, Michigan, US, which includes most of Detroit's metropolitan area. As in \cite{heleno2022optimizing}, the tracts considered in this case study qualify to the US Department of Energy (DOE) Communities Local Energy Action Planning (LEAP) program \cite{U.S.DepartmentofEnergy2021}, which requires a tract to have a low-income population $\geq$ 30\% and median household energy burden $\geq$ 6\%. Furthermore, the tract must meet at least one of the following criteria: 1) the community has an historical economic dependence on fossil fuel industrial facilities including extraction, processing, or refining; or 2) the tract is classified as experiencing moderate or high susceptibility on the U.S. Environmental Protection Agency’s Environmental Justice Screening (EJSCREEN) tool \cite{U.S.EnvironmentalProtectionAgency2021}.

 Similarly to \cite{heleno2022optimizing}, the set of households considered meet the criteria of having an annual income ranging from 80\% to 100\% of the Area Median Income (AMI). However, since we are validating our model against a full time-resolved formulation, we had to reduce the size of the household archetypes dataset, to account for the computational burden of the benchmark solution. Therefore, in comparison with the case study used in \cite{heleno2022optimizing}, only rental households are included, resulting in a total of 3,651 households represented by 1,019 household archetypes across 436 census tracts.  
 
The cost and economic assumptions, including battery costs and lifetimes, are presented in Table~\ref{tab:intervention_costs_case_study}. The simulation assumes a discount rate of 3\% and a net-billing scheme in which the PV sales are compensated at 60\% of the electricity demand rate ($PV_{rem} = 0.6 \cdot P_{el}$).

\begin{table}[!ht]
    \centering
    \begin{tabular}{ccc}
    \hline
    Intervention & Cost (/MW or /building) & Lifetime \\
    \hline
    Rooftop Solar &	  \$2.4M & 20 \\
   Community Solar &  \$1.6M & 20	  \\
   Community Wind &   \$2.5M&	 15  \\
   Battery (4-hour) & \$1.2M& 5 \\
   Weatherization & Same costs per building type as in \cite{heleno2022optimizing} &	 35\\
    \hline
    \end{tabular}
    \caption{Cost of interventions considered in the analysis.}
    \label{tab:intervention_costs_case_study}
\end{table}

\subsection{Methodology validation}

To validate the methodology proposed in this paper, which introduces a linear approximation of the temporal representation of storage operations, we compare it with a fully time-dependent model, which incorporates hourly load and photovoltaic profiles for each household archetype and considers the full set of equations described in (\ref{eq:stagg})-(\ref{eq:lim_dc}). For simplicity, in this result section, we refer to the proposed approach as the "LINEARIZED" model, whereas the fully time-dependent model is referred to as the "TIME" model. It is important to note that the only difference between these models is the time-resolved operational constraints, as in both cases the objective is to provide the optimal portfolio of energy interventions that minimize both inequity across the targeted households and the budget required for implementation.

We run both models considering the set of 1,019 household archetypes with an unconstrained budget. From a computational time perspective, the performance of the models cannot be compared. As shown in (\ref{eq:stagg})-(\ref{eq:lim_dc}, the fully time-resolved "TIME" model requires state-of-charge, charging, and discharging variables for each of the 8,760 hours of the year per archetype. When running the optimization for the 1,019 household archetypes, the model took several hours to converge. In contrast, the proposed approximation, which includes only one new binary variable per archetype, solved in a few seconds.

In terms of accuracy, Table \ref{tab:table_T} compares the resulting values of the inequity objective function and the consequent budget, along with the average energy burden and average energy burden reduction. The errors for all equity and economic variables are less than 5\%, showing that the approximate representation of the energy operation leads to a very similar impact when compared with a fully temporally resolved model.

\begin{table}[h!]
    \centering
    \resizebox{\textwidth}{!}{%
    \begin{tabular}{cccc}
        \hline
        \textbf{Indicator} & \textbf{Time Model} & \textbf{Linearized Model} & \textbf{Error} \textbf{(in \%)}\\
        \hline
        Inequity Function (\%)& 36483.4 & 34741.3 & 4.8\\ 
        Budget Function (M\$) &  2.748 &  2.671 & 2.8\\
        Average Energy Burden After Interventions (\%) & 5.11 & 5.16 & 1.0\\ 
        Average $\Delta eb_k^+$ & 0.099 & 0.0951 & 3.9\\  
        \hline
    \end{tabular}}
    \caption{Key comparison indicators between the time-dependent and the linearized model}
    \label{tab:table_T}
\end{table}

\begin{figure}
    \centering
    \includegraphics[scale = 0.4]{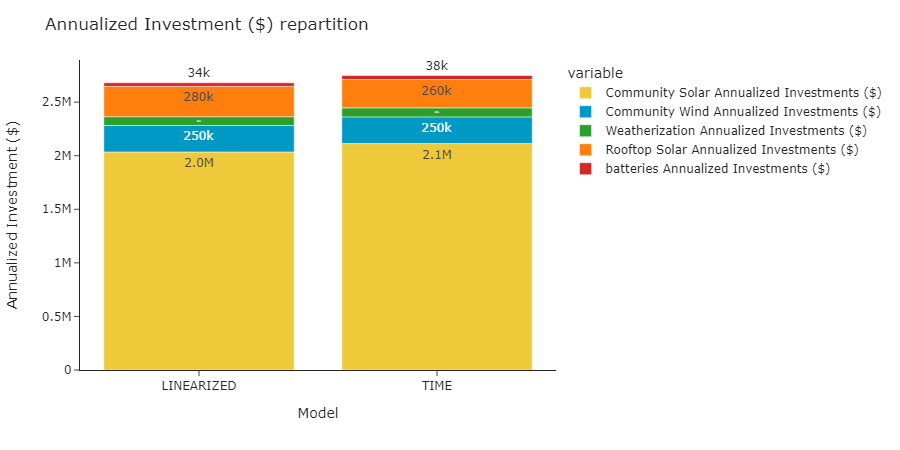}
    \caption{Annualized investments corresponding to the optimal portfolio of energy interventions for the linearized and the full time dependent models, considering $PV_{rem} = 0.6 \cdot P_{el}$.}
    \label{fig:comp_costs}
\end{figure}

\begin{figure}
        \centering
        \includegraphics[scale = 0.4]{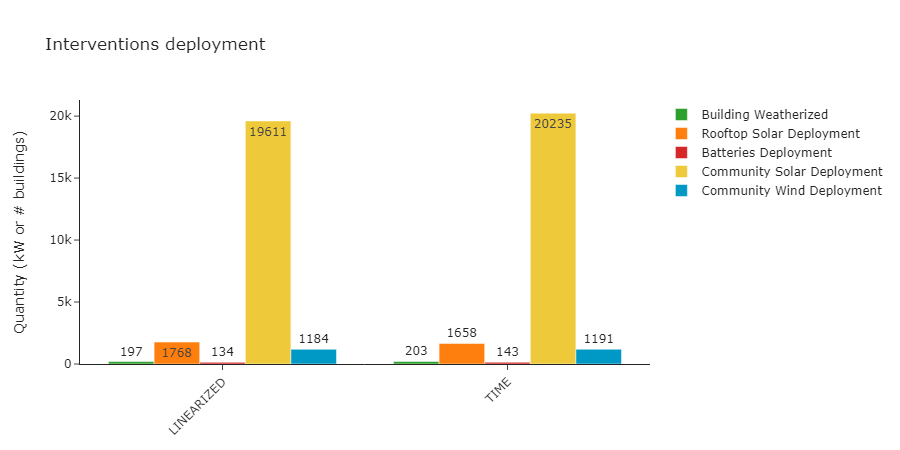}
        \caption{Optimal portfolio of energy interventions for the linearized and the full time dependent models, considering $PV_{rem} = 0.6 \cdot P_{el}$.}
        \label{fig:comp_solutions}
    \end{figure}

Figure~\ref{fig:comp_costs} compares the resulting optimal portfolio of investments that minimize the inequity objective function in both models. Figure~\ref{fig:comp_solutions} provides an analogous comparison for the corresponding intervention types. As shown in both figures, the two models offer very similar optimal portfolio solutions for equitable investments, demonstrating the accuracy of the approximated linearized solution presented in this paper. However, it is possible to find errors in the linearized model solution associated with the approximation of the storage operation. For example, in Figure \ref{fig:comp_solutions}, it is interesting to note a difference in the ratio between rooftop PV and storage capacities. While the fully time-resolved model finds an optimal solution with 1,658 kW of rooftop PV and 143 kW of batteries (86 kW of storage per MW of PV), the linearized model finds a solution with 1,768 kW of rooftop PV and 134 kW of batteries (76 kW of storage per MW of PV). This result suggests that the linearized model may be overestimating the amount of solar that can be stored in a battery due to the linear approximation previously illustrated in Figure~\ref{figure:solar_sc_linearization}.

\begin{figure}
        \centering
        \includegraphics[scale = 0.43]{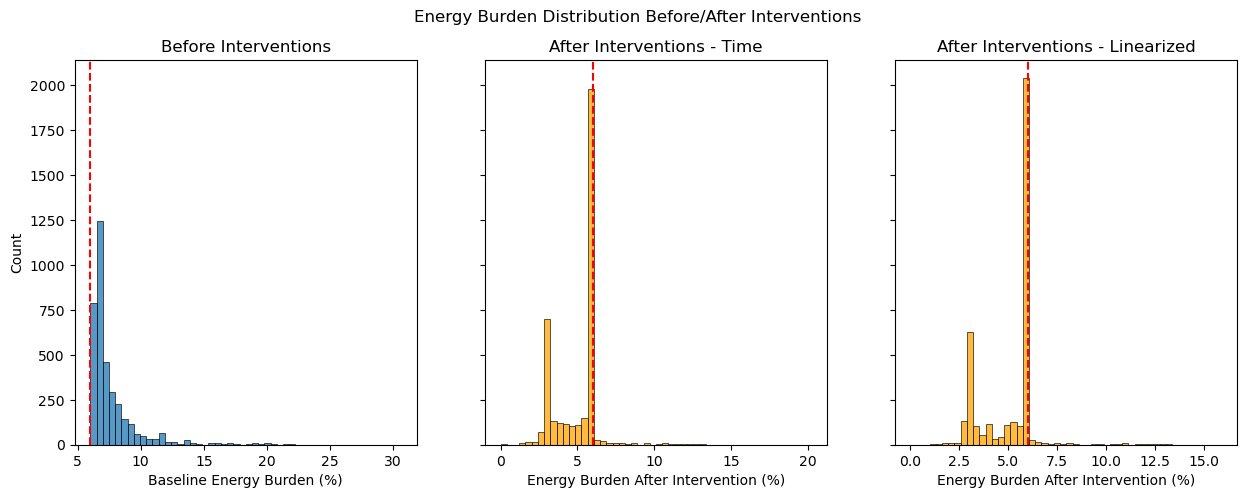}
        \caption{Energy Burden Before and after the interventions for the linearized and the full time dependent models.}
        \label{fig:Eb_ba_T}
\end{figure}

As a final step in validation, we compare how different battery modeling approaches impact the energy burden reduction of individual household archetypes. Figure~\ref{fig:Eb_ba_T} plots the distribution of the energy burden of the targeted population before and after the deployment of weatherization and DER interventions. As expected, in the left panel, it is possible to observe that before the interventions, the entire population had an energy burden above the 6\% threshold, with some household archetypes reaching extreme burden conditions (around 20\%). After the interventions, achieved with both the linearized and the full time-resolved models, almost the entire population falls under the 6\% energy burden threshold. When comparing the resulting energy burden distributions, no major changes can be observed between the linearized and the full temporal model. This result confirms that the methodology proposed to represent storage operations in equitable planning can accurately approximate the full time-resolved model.

\subsection{Storage deployment}

\begin{figure}[h!]
    \centering
    \includegraphics[scale = 0.4]{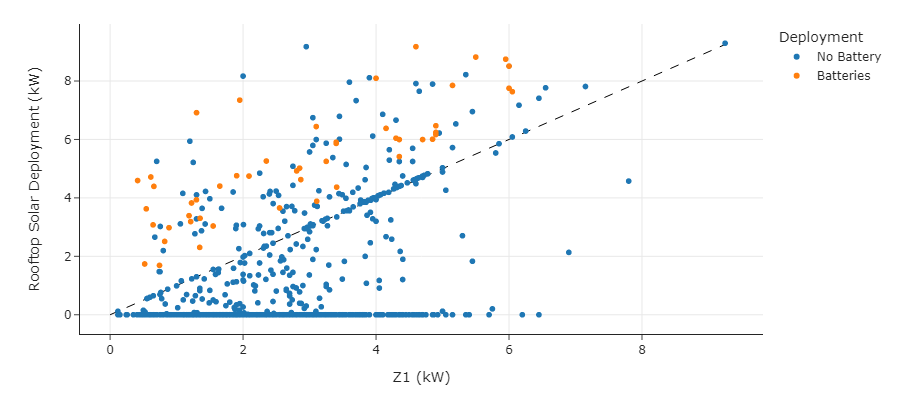}
    \caption{Rooftop solar and corresponding storage deployment as a function of the threshold variable Z1, considering $PV_{rem}=0.6 \cdot Pel$.}
    \label{fig:batt_dep}
\end{figure}

In the model presented above, storage is deployed to maximize the revenue of rooftop solar assets to help reduce the overall energy burden of a household. Hence, in this subsection, we use our linearized model to understand how storage is deployed as a function of the rooftop PV in households to support equity objectives. Figure~\ref{fig:batt_dep} presents, for each household archetype, the rooftop solar intervention as a function of the threshold variable Z1, introduced in section \ref{sec:methodology}, that indicates the size above which the PV array starts generating surplus. Thus, rooftop solar installations plotted under the diagonal have a small capacity relative to the load and will be used for self-consumption only, whereas the ones above the diagonal will be large enough to generate surplus. As shown in the Figure, batteries are only deployed with the latter group, since without surplus there is no incentive to store energy. 

\begin{figure}
    \centering
   \begin{subfigure}{0.5 \textwidth}
        \centering
        \includegraphics[width=0.95\linewidth]{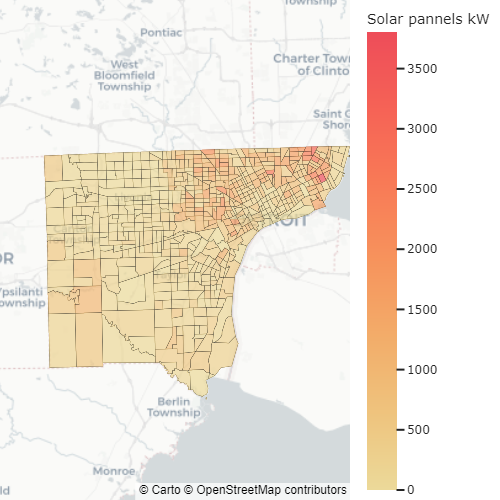}
        \caption{Distribution of rooftop interventions.}
        \label{fig:solar_locational}
   \end{subfigure}%
   \begin{subfigure}{0.5\textwidth}
        \centering
        \includegraphics[width=0.95\linewidth]{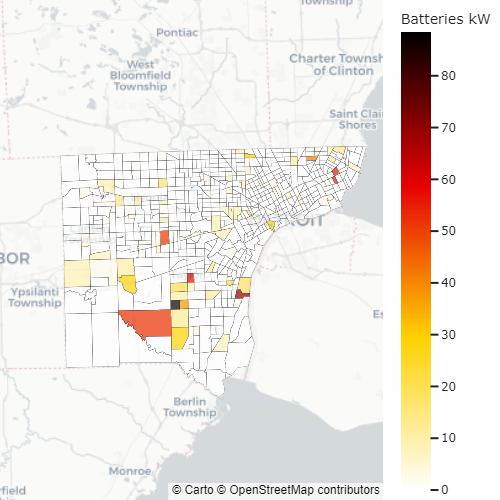}
        \caption{Distribution of storage interventions.}
        \label{fig:storage_locational}
   \end{subfigure}  
   \caption{Deployment of rooftop solar and storage interventions for different census tracts.}
    \label{fig:battery_solar_comparision}
\end{figure}

In Figure~\ref{fig:batt_dep}, it is also possible to observe that not all rooftop PV installations with surplus (above the diagonal) are paired with a battery, since the deployment of "solar only" solutions may be enough to address the energy burden of certain households. In fact, this becomes evident with the comparison of batteries and rooftop PV deployment across the territory, as shown in Figure \ref{fig:battery_solar_comparision}. While rooftop PV is deployed throughout the county, batteries are only deployed in certain census tracts. In addition, it is important to emphasize that the tracts with more storage deployment are not necessarily the ones with more rooftop PV. This shows that, unlike rooftop solar, batteries get deployed strategically to reduce specific situations of extreme energy burden not as a general solution for the problem. The reason is that solar PV (both in its rooftop or community form) is a low cost solution to directly offset load consumption. In contrast, storage in an expensive asset for the value it creates, which is solely the ability to transfer PV surplus from the moments when it is less valuable to later periods where it can be used to offset home consumption.

\begin{figure}
    \centering
    \includegraphics[scale = 0.6]{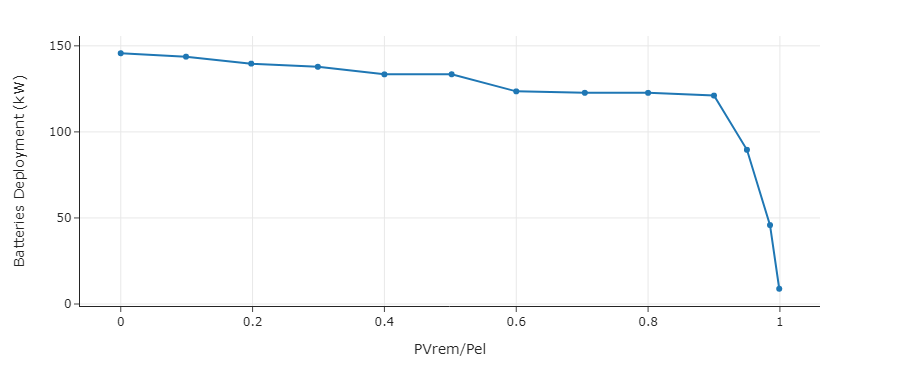}
    \caption{Storage deployment as a function of the solar compensation represented by the $PV_{rem}/Pel$ ratio.}
    \label{fig:sens_batt}
\end{figure}

To better illustrate this arbitrage value of storage associated with the solar compensation scheme, we compared the deployment of storage for different ratios between PV sales ($PV_{rem}$) and electricity demand rates (${Pel}$). Figure~\ref{fig:sens_batt} shows that as soon as there is a net-billing scheme ($PV_{rem}<Pel$), there is an opportunity value of shifting PV surplus in time, and therefore storage gets deployed. When the solar is compensated in a net-metering scheme ($PV_{rem}=Pel$), there is no opportunity for storage to decrease the energy burden at the household level.

\section{Conclusion}
\label{sec:conclusion}

This work has successfully extended the Justice40 model to incorporate energy storage systems in the portfolio of place-based equity investments. We achieve this by proposing a linear approximation of the temporal representation of storage operations, which can significantly reduce computational complexity while maintaining high accuracy in results. This approach was validated against a full time-resolved model and the comparison of equity and economic results showed similar optimal investment portfolio solutions and energy energy burden reduction impacts. 

As expected, storage becomes a viable solution to decrease energy burden under net-billing schemes where PV surplus can be shifted to more valuable periods. However, our results highlight that  batteries are allocated more selectively to address specific cases of extreme energy burden, unlike rooftop PV, which is deployed more broadly as an universal solutions.

Potential future works could extend the storage modeling approach proposed in this paper beyond solar compensation schemes and explore other equitable valuates of energy storage, such as resilience. This can potentially allow, for example, to impose a resilience constraints in the selection of the investments in the Justice40 model.

\bibliography{references}

\end{document}